\documentclass[12pt]{amsart}

\usepackage{amsmath,amsfonts}
\usepackage{amsthm}
\usepackage{enumerate}

\newtheorem{theorem}{Theorem}[section]
\newtheorem{proposition}[theorem]{Proposition}
\newtheorem{lemma}[theorem]{Lemma}
\newtheorem{corollary}[theorem]{Corollary}

\theoremstyle{definition}
\newtheorem{definition}[theorem]{Definition}
\newtheorem{example}[theorem]{Example}

\theoremstyle{remark}
\newtheorem{remark}[theorem]{Remark}

\numberwithin{equation}{section}

\newcommand{\nC}{\mathbb C}

\newcommand{\nN}{\mathbb N}
\newcommand{\nR}{\mathbb R}

\newcommand{\cD}{{\mathcal D}}

\newcommand{\cK}{{\mathcal K}}
\newcommand{\cL}{{\mathcal L}}

\newcommand{\cP}{{\mathcal P}}

\newcommand{\cV}{{\mathcal V}}
\newcommand{\cU}{{\mathcal U}}
\newcommand{\cZ}{{\mathcal Z}}

\newcommand{\mb}[1]{\mathbf{#1}}

\newcommand{\dist}{d} 
 
\newcommand{\ds}{d_F}
\newcommand{\dpo}{d_{\cP}} 

\begingroup
\def\2#1{\ifnum#1<10 0\fi\the#1}
\xdef\isodayandtime%
{\the\year-\2\month-\2\day\space\2{\count0}:\2{\count2}}

\endgroup

\begin{document}

\title[Roots and Polynomials as Homeomorphic Spaces]%
{Roots and Polynomials \\
      as Homeomorphic Spaces}

\author{Branko \'{C}urgus}
\address{Department of Mathematics, Western Washington
University, \newline \hspace*{3mm} Bellingham, WA 98225}
\email{curgus@cc.wwu.edu}
\author{Vania Mascioni}
\address{Department of Mathematical Sciences, Ball State
University,\newline \hspace*{3mm} Muncie, IN 47306-0490}
\email{vdm@cs.bsu.edu}


\date{October 30, 2004}


\begin{abstract}
\hspace*{-2.1mm} We provide a unified, elementary, topological
approach to the classical results stating the continuity of the
complex roots of a polynomial with respect to its coefficients,
and the continuity of the coefficients with respect to the roots.
In fact, endowing the space of monic polynomials of a fixed degree
$n$ and the space of $n$ roots with suitable topologies, we are
able to formulate the classical theorems in the form of a
homeomorphism. Related topological facts are also considered.
\end{abstract}

\maketitle

The roots of a polynomial depend continuously on its coefficients.
This is probably the best known and most used perturbation theorem
and, clearly, it is a continuity statement (see \cite{M} for
several historical references; also, see our final remarks in
Section \ref{sfr}). Conversely, the coefficients depend
continuously on the roots. This is essentially due to Vi\`ete's
formulas; see Theorem \ref{tvpc} below. However, this second
result is often formulated separately from the first, and there
has been no unanimity as to the topology on the set of roots.

In this note we provide a metric space setting in which both of
these classical continuity results can be stated as a unique
homeomorphism (our main result will be Theorem \ref{tmt}) between
the corresponding metric spaces of roots and polynomials. This
reveals more than may be widely known about the similar
topological structure of these spaces.

We only use very basic background knowledge of the topology of
metric spaces for example at the level of Rudin's  or Baum's
classical books \cite{R,B}.  Whenever we refer to a set as a
metric space we imply that a specific metric has been earlier
defined on it. Each subset of a metric space is considered a
metric space with the induced metric.  We use the standard
notation $\nN$ for the set of positive integers, $\nR$ for the set
of real numbers, $\nC$ for the set of complex numbers,  and $i =
\sqrt{-1}$. Throughout this note $n\geq 2$ is a fixed positive
integer. We study complex monic polynomials of order $n$ and we
consider all their complex roots. Since monic polynomials of
degree one are in an obvious one-to-one correspondence with their
unique root, the case $n=1$ is a special, though trivial, case.
Note that Theorem \ref{tKnc} and Corollary \ref{cnc} are not true
in the case $n=1$.

\section{Metric Space Preliminaries} \label{smsp}

\begin{definition}
Let $(X,d_X)$ and $(Y,d_Y)$ be metric spaces, and let $f:X\to Y$
be a bijection. If both $f$ and $f^{-1}$ are continuous then $f$
is called a {\it homeomorphism} between $X$ and $Y$.
\end{definition}

Our first theorem bears a strong resemblance to the classical
result that states that \textit{a continuous bijection from a
compact space to a Hausdorff space has a continuous inverse} (see
\cite[Theorem 3.21]{B} or \cite[Theorem 4.17]{R}, for example).

%
%

\begin{theorem} \label{tc}
Let  $(X,d_X)$ and $(Y,d_Y)$  be metric spaces and let $f : X
\rightarrow Y$ be a bijection. Suppose that the following three
conditions are satisfied:
\begin{enumerate}[{\rm (a)}]
\item \label{tca}
Each bounded and closed subset of $X$ is compact.
\item \label{tcb}
$f$ is continuous.
\item \label{tcc}
$f^{-1}$ maps each bounded set in $Y$ into a bounded set in $X$.
\end{enumerate}
Then $f^{-1}$ is continuous (and so $f$ is a homeomorphism).
\end{theorem}
\begin{proof}
Let $\{y_k\}$ be a convergent sequence in $(Y,d_Y)$ with limit
$y$. Since $\{y_k\}$ is bounded, assumption (\ref{tcc}) implies
that the sequence $\bigl\{f^{-1}(y_k)\bigr\}$ is bounded in $X$
and thus it is contained in a closed and bounded subset of $X$.
Recall that in a metric space if a set is bounded, that is if it
has a finite diameter, then its closure has the same diameter. By
(\ref{tca}), $\bigl\{f^{-1}(y_k)\bigr\}$ has a convergent
subsequence. If $\bigl\{f^{-1}(y_{m_k})\bigr\}$ is an arbitrary
convergent subsequence of $\bigl\{f^{-1}(y_k)\bigr\}$ with, say,
\begin{equation*}
\lim\limits_{k\rightarrow \infty} f^{-1}(y_{m_k}) = x
\end{equation*}
the continuity of $f$ (assumption (\ref{tcb})) implies that
\begin{equation*}
\lim\limits_{k\rightarrow \infty} y_{m_k} = f(x) =
\lim\limits_{k\rightarrow \infty} y_{k}  = y.
\end{equation*}
Thus each convergent subsequence of the bounded sequence
$\bigl\{f^{-1}(y_k)\bigr\}$ converges to the same element
$f^{-1}(y)$, and this implies that $\bigl\{f^{-1}(y_k)\bigr\}$
converges to $f^{-1}(y)$. Since the sequence $\{y_k\}$ was an
arbitrary convergent sequence in $Y$, the theorem is proved.
\end{proof}

\begin{proposition} \label{p13}
If each bounded and closed subset of a metric space $(X,d_X)$ is
compact, then $(X,d_X)$ is complete.
\end{proposition}
\begin{proof}
Each Cauchy sequence in a metric space is bounded and thus contained
in a closed ball. Since by assumption a closed ball in $(X,d_X)$ is
compact, each Cauchy sequence in $(X,d_X)$ has a convergent
subsequence. Consequently each Cauchy sequence in $(X,d_X)$
converges.
\end{proof}

By $\nC^n$ we denote the set of all ordered $n$-tuples of complex
numbers. We equip this space with what is called the ``supremum
norm''
\begin{equation*}
\|{\mathbf v}\|_{\infty} = \max_{1\leq j\leq n} |v_j|
 \ \ \ \ \text{for} \ \ \ \
{\mathbf v} = \bigl(v_1,\ldots,v_n\bigr) \in \nC^n,
\end{equation*}
and, for ${\mathbf u}, {\mathbf v} \in \nC^n$, ${\mathbf u} =
(u_1,\ldots,u_n), \  {\mathbf v} = (v_1,\ldots,v_n)$, the
corresponding metric
\begin{equation*}
d_{\infty}({\mathbf u},{\mathbf v}) =
 \max_{1\leq j \leq n} |u_j-v_j|
  = \|{\mathbf u}-{\mathbf v} \|_{\infty} .
\end{equation*}
The following proposition is well known and not difficult to
prove.
\begin{proposition}[Heine-Borel] \label{pCn}
$(\nC^n,d_{\infty})$ is a metric space. A subset of
$(\nC^n,d_{\infty})$ is compact if and only if it is bounded and
closed.
\end{proposition}

The metric $d_{\infty}$ on $\nC^n$ is chosen for convenience only.
Clearly it can be replaced with any other equivalent metric.

Next we prove a topological property of the space
$(\nC^n,d_{\infty})$ which we shall need in Section \ref{srCn}.

\begin{definition} \label{dpwc}
Let $(X,d_X)$ be a topological space. A subset $S$ of $X$ is {\em
pathwise connected} if for each $u,v \in S$ there exists a
continuous function $\Theta:[0,1] \to X$ such that $\Theta(0) = u,
\Theta(1) = v$, and the range of $\Theta$ is a subset of $S$. The
range of the function $\Theta$ is called a {\em path} from $u$ to
$v$ which is contained in $S$.
\end{definition}

\begin{lemma} \label{lDpwc}
Let $\cD$ be the subset of \ $\nC^n$ consisting of all $n$-tuples
of distinct complex numbers. Then $\cD$ is an open
pathwise-connected subset of $(\nC^n,d_{\infty})$.
\end{lemma}
\begin{proof}
Given the continuous function $f(z_1,...,z_n) = \prod_{i\not=j}
(z_i-z_j)$ between $\nC^n$ and $\nC$, we can write
$\cD=f^{-1}(\nC\!\setminus\!\{0\})$, and since
$\nC\!\setminus\!\{0\}$ is open in $\nC$, $\cD$ must be open in
$\nC^n$.

To prove that $\cD$ is pathwise connected, we let
$\mb{v}=(v_1,\ldots,v_n)$ and $\mb{w}=(w_1,\ldots,w_n)$ be two
points in $\cD$ and construct a path from $\mb{v}$ to $\mb{w}$ which
is contained in $\cD$.

First consider a special case. Assume that there exists $k \in
\{1,\ldots,n\}$ such that $v_j = w_j$ for all $j \in
\{1,\ldots,n\}\!\setminus\! \{k\}$ and $v_k \neq w_k$.  Since the
numbers $v_k, w_k, v_1, \ldots, v_{k-1}, v_{k+1}, \ldots, v_{n}$
are mutually distinct, it is not hard to construct a continuous
function $\phi:[0,1] \to \nC$ such that $\phi(0) = v_k,\, \phi(1)
= w_k$ and none of the numbers $v_1, \ldots, v_{k-1}, v_{k+1},
\ldots, v_{n}$ is in the range of $\phi$.  Consequently the
function
 \[
\Theta(t) = \bigl( v_1, \ldots, v_{k-1}, \phi(t) , v_{k+1}, \ldots,
v_{n} \bigr), \ \ \ \ t \in [0,1],
 \]
is a path from $\mb{v}$ to $\mb{w}$ which is contained in $\cD$.

Now consider the general case of arbitrary points
$\mb{v}=(v_1,\ldots,v_n)$ and $\mb{w}=(w_1,\ldots,w_n)$ in $\cD$.
Let $\mb{u}=(u_1,\ldots,u_n) \in \cD$ be such that
\begin{equation*} \label{eqvuw}
\bigl\{ u_1,\ldots,u_n \bigr\} \, \cap \, \bigl\{ v_1,\ldots,v_n,
w_1,\ldots,w_n \bigr\} \, = \, \emptyset.
\end{equation*}
Consider the following sequence of points in $\cD$:
 \begin{align*}
 & (v_1,v_2,v_3,\ldots,v_{n-1},v_n), \
      (u_1,v_2,v_3,\ldots,v_{n-1},v_n), \\
 & (u_1,u_2,v_3,\ldots,v_{n-1},v_n), \ \ldots, \
         (u_1,u_2,u_3,\ldots,u_{n-1},v_n), \\
 & (u_1,u_2,u_3,\ldots,u_{n-1},u_n), \
            (w_1,u_2,u_3,\ldots,u_{n-1},u_n), \\
 & (w_1,w_2,u_3,\ldots,u_{n-1},u_n), \ \ldots, \
              (w_1,w_2,w_3,\ldots,u_{n-1},u_n), \\
 & (w_1,w_2,w_3,\ldots,w_{n-1},u_n), \ (w_1,w_2,w_3,\ldots,w_{n-1},w_n).
 \end{align*}
The special case considered above applies to each of the $2n$ pairs
of consecutive points in this sequence.  It follows that for each of
these pairs there exists a path contained in $\cD$ which connects
them.  Since each two consecutive pairs contain a point in common,
these $2n$ paths connect to a path connecting $\mb{v}$ and $\mb{w}$
which is clearly contained in $\cD$.  As $\mb{v}$ and $\mb{w}$ were
arbitrary points in $\cD$ this proves that $\cD$ is pathwise
connected.
\end{proof}

By $\cP_{n,1}$ we denote the set of all monic complex polynomials
of degree $n$.  Let
\begin{equation*}
f(z) =  z^n + a_{n-1} z^{n-1} + \cdots + a_0 , \ \ g(z) = z^n +
b_{n-1} z^{n-1} + \cdots + b_0, \ \ z \in \nC,
\end{equation*}
be in $\cP_{n,1}$.  Define a metric on $\cP_{n,1}$ by
\begin{equation} \label{poltop}
\dpo(f,g) : = \max\bigl\{ |a_0 - b_0|, \ldots, |a_{n-1} - b_{n-1}|
\bigr\} \,.
\end{equation}

\begin{proposition} \label{ePn}
$(\cP_{n,1},\dpo)$ is a metric space. A subset of the metric space
$(\cP_{n,1},\dpo)$ is compact if and only if it is bounded and
closed.
\end{proposition}
\begin{proof}
The function
\begin{align*}
\bigl(v_1,\ldots,v_{n} \bigr) \ \longmapsto \ p \ \ \
 & \text{where} \ \ \ \ p(z) = z^n + v_{n} z^{n-1} + \cdots + v_1, \\
 & \ \ \ \text{and} \ \ \ \ \bigl(v_1,\ldots,v_{n} \bigr) \in
 \nC^n,
\end{align*}
is a distance preserving bijection between the spaces
$(\nC^n,\dist_{\infty})$ and $(\cP_{n,1},\dpo)$. Therefore the
proposition follows from Proposition \ref{pCn}.
\end{proof}

\section{The Metric Space of Roots} \label{smsor}

At the end of Section \ref{smsp} we introduced the metric space
$(\cP_{n,1},\dpo)$ of all monic polynomials of degree $n$.  Now we
define the space of sets of roots of these polynomials.  Since
roots can occur with finite multiplicities, instead of the set of
roots of a polynomial we consider the multiset of roots, that is,
we allow elements to occur with multiplicities.  Denote by $\cZ_n$
the family of all multisets of complex numbers with $n$ elements.
For multisets $U = \{u_1,\dotsc,u_n\}$ and $V =
\{v_1,\dotsc,v_n\}$ in $\cZ_n$, define
\begin{equation} \label{eqdF}
\ds(U,V) := \min_{\tau\in\Pi_n}  \max_{1\leq j\leq n} |u_j -
v_{\tau(j)}|,
\end{equation}
where $\Pi_n$ is the set of all permutations of $\{1,\dotsc,n\}$.
The function $\ds$, which is a metric by the proposition below, is
analogous to the Fr\'{e}chet metric defined for curves in
\cite[Chapter 6]{E}.  Instead of curves here we have multisets and
a function
 $
f :\{1,\ldots,n \} \to \nC
 $
is a parametrization of the multiset $\{f(k) : 1 \leq k \leq n\}$.
If we denote by $\cU$ and $\cV$ all possible parameterizations of
multisets $U$ and $V$, respectively, then definition \eqref{eqdF}
can be rewritten as
\begin{equation*} 
\ds(U,V) = \min_{f\in\cU,g\in\cV}  \max_{1\leq k\leq n} |f(k) -
g(k)|.
\end{equation*}

\begin{proposition} \label{Ldistine}
The function \ $\ds: \cZ_n \times \cZ_n \rightarrow [0,\infty)$ is
a metric on $\cZ_n$.
\end{proposition}
\begin{proof}
Let $U, V, W \in \cZ_n$. We need to prove the following three
properties of $\ds$:
\begin{gather} \label{pr0}
\ds(U,V) = 0 \ \ \Longleftrightarrow \ \ U = V,
 \\ \label{sym} \ds(U,V) = \ds(V,U),
 \\ \label{trs}%
\ds(U,V) \leq \ds(U,W) + \ds(W,V) \, .
\end{gather}
To prove \eqref{pr0} is a simple exercise. The definition of $\ds$
can be rewritten as
\begin{equation} \label{eqdss}
\ds(U,V) = \min_{\sigma,\tau\in\Pi_n} \max_{1\leq j\leq n}
|u_{\sigma(j)} - v_{\tau(j)}| .
\end{equation}
Since the last expression is symmetric in $U$ and $V$, this shows
that $\ds(U,V) = \ds(V,U)$ and thus \eqref{sym} holds.

To prove \eqref{trs} note that the triangle inequality for complex
numbers yields
\begin{equation} \label{distine3}
|u_j - v_{\tau(j)}| \leq |u_j - w_{\sigma(j)}| + |w_{\sigma(j)} -
v_{\tau(j)}|,
\end{equation}
for arbitrary $j \in \{1,\dotsc,n\}$ and arbitrary $\sigma, \tau
\in \Pi_n$. Keeping $\sigma$ and $\tau$ fixed and taking maximums
with respect to $j \in \{1,\ldots,n\}$ in \eqref{distine3} we get
\begin{equation} \label{distine4}
\max_{1\leq j\leq n} |u_j - v_{\tau(j)}| \leq \max_{1\leq l \leq
n} |u_l - w_{\sigma(l)}| + \max_{1\leq k \leq n} |w_{\sigma(k)} -
v_{\tau(k)}|.
\end{equation}
Keeping $\sigma \in \Pi_n$ fixed and taking the minimums of both
sides in \eqref{distine4} with respect to $\tau \in \Pi_n$ we get
\begin{equation*} 
\ds(U,V) \leq \max_{1\leq l \leq n} |u_l - w_{\sigma(l)}| + \ds(W,V)
\end{equation*}
and so \eqref{trs} follows by taking the minimum of the right hand
side with respect to $\sigma \in \Pi_n$.
\end{proof}

Next we explore the relationship between the space $(\cZ_{n},\ds)$
and the more familiar space $(\nC^n,\dist_{\infty})$. First we
define two functions $P$ and $K$.

Define $P: \nC^n \rightarrow \cZ_n$ by
\begin{equation} \label{eqdP}
P \bigl(\, ( v_1,\ldots,v_n ) \, \bigr) := \bigl\{
v_1,\ldots,v_n\bigr\}, \ \ \ \ ( v_1,\ldots,v_n ) \in \nC^n.
\end{equation}
Here an $n$-tuple is simply mapped to the multiset of its elements
(once again, with multiplicities preserved). By the definitions of
$\ds$ and $\dist_{\infty}$ it follows that
\begin{equation} \label{dsineq}
\ds\bigl(P({\mathbf v}),P({\mathbf w}) \bigr) \leq \|{\mathbf
v}-{\mathbf w}\|_{\infty} =  \dist_{\infty}({\mathbf v},{\mathbf
w}) \ \ \ \text{for all} \ \ \ {\mathbf v},{\mathbf w} \in \nC^n.
\end{equation}
Thus $P: \nC^n \rightarrow \cZ_n$ is a contraction (and therefore
a continuous function) between $(\nC^n,\dist_{\infty})$ and
$(\cZ_n,\ds)$.

Clearly $P$ is onto, but not one-to-one. For each $V \in \cZ_n$
the set
 \[
P^{-1}(V) : = \{{\mathbf v} \in \nC^n : P({\mathbf v}) = V\}
 \]
has between $1$ and $n!$ elements, depending on the multiplicities
of the elements in $V$.  Note that for distinct $V$ and $W$ in
$\cZ_n$ the sets $P^{-1}(V)$ and $P^{-1}(W)$ are disjoint.

To define a partial inverse of $P$ let $\cK$ be a subset of $\nC^n$
with the property that for each $V \in \cZ_n$ the set $\cK \cap
P^{-1}(V)$ has exactly one element. (In Example \ref{elo1} below we
give a specific example of a set $\cK$ with this property.) This
assumption is equivalent to the requirement that the restriction
\begin{equation*}
P|_{\cK}: \cK \to \cZ_n
\end{equation*}
of $P$ onto $\cK$ is a bijection. In this way to each $V =
\{v_1,\ldots,v_n\} \in \cZ_n$ we associate a unique $n$-tuple
$(v_1,\ldots,v_n) \in \nC^n$ that has exactly the elements of $V$
as coordinates.  Now define the function $K: \cZ_n \to \nC^n$ by
 \begin{equation} \label{eqdK}
K : = \bigl( P|_{\cK}\bigr)^{-1}.
 \end{equation}
As an immediate consequence of the definitions we conclude that
$P\circ K$ is the identity on $\cZ_n$.

Let $O \in \cZ_n$ be the multiset consisting of $n$ zeros. By the
definitions of $\ds$ and $K$ it follows that
\begin{equation} \label{zeroset}
 \ds\bigl(V,O\bigr) = \|K(V)\|_{\infty} \ \ \ \
\text{for all} \ \ \ V \in \cZ_n.
\end{equation}

\begin{proposition} \label{pbdds}
Let $\{V_k\}$ be a sequence in $(\cZ_n,\ds)$. The following
statements are equivalent.
\begin{enumerate}[{\rm (a)}]
\item The sequence $\{V_k\}$ is bounded in $(\cZ_n,\ds)$. \item
The set $\bigcup_{k=1}^{\infty} V_k$ of complex numbers is bounded
in $\nC$.
\item The sequence $\{K(V_k)\}$ is bounded in
$(\nC^n,\dist_{\infty})$.
\end{enumerate}
\end{proposition}
\begin{proof}
Let $\{V_k\}$ be a bounded sequence in $(\cZ_n,\ds)$. Since
$\{V_k\}$ is bounded  there exists $M > 0$ such that
\begin{equation}
\ds\bigl(O,V_k \bigr) < M \ \ \ \ \text{for all} \ \ k \in \nN.
\end{equation}
By \eqref{zeroset} part (b) follows trivially, and just as
trivially (b) implies (c). If (c) holds, then \eqref{zeroset}
implies that the sequence $\{\ds(V_k,O)\}$ is bounded, and thus
(a) follows just as easily.
\end{proof}
In a similar way \eqref{zeroset} can be used to prove the
following proposition.
\begin{proposition} \label{pKbb}
The function $K: \cZ_n \to \nC^n$ maps each bounded set in
$(\cZ_n,\ds)$ to a bounded set in $(\nC^n,\dist_{\infty})$.
\end{proposition}
The continuity of $K$ is discussed in Section \ref{srCn} (see, in
particular, Corollary \ref{cnc}).
\begin{theorem} \label{tbscs}
A subset of the metric space $(\cZ_n,\ds)$ is compact if and only
if it is bounded and closed.
\end{theorem}
\begin{proof}
Let $\cV$ be an arbitrary bounded and closed subset of $\cZ_n$. To
prove that $\cV$ is compact we shall prove that an arbitrary
sequence $\{V_k\}$ in $\cV$ has a convergent subsequence.  By
Proposition \ref{pbdds} the sequence $\{K(V_k)\}$ is bounded in
$(\nC^n,\dist_{\infty})$.  By the Bolzano-Weierstrass Theorem
there exists a subsequence $\{V_{m_k}\}$ of $\{V_k\}$ such that
$\{K(V_{m_k})\}$ converges, say, to the $n$-tuple ${\mathbf w} =
(w_1,\ldots,w_n)$, in $(\nC^n,\dist_{\infty})$. Since $P: \nC^n
\rightarrow \cZ_n$ is continuous and $P\circ K$ is the identity on
$\cZ_n$, it follows that $\{ V_{m_k}\}$ converges to $P({\mathbf
w})$ in $(\cZ_n,\ds)$. Since $\cV$ is closed $P({\mathbf w}) \in
\cV$, and thus $\cV$ is compact. Since the converse is true in
each metric space the theorem is proved.
\end{proof}

In the next two examples we use the notion of lexicographic
ordering in $\nC$.  Let $a,b,c,d \in \nR$. For two complex numbers
$a+ib$ and $c+id$ the lexicographic ordering $a+ib \preceq c+id$
is defined by
\begin{equation*}
a+ib \preceq c+id \ \ \Longleftrightarrow \ \ \bigl[\bigl(a <
c\bigr) \ \ \vee \ \ \bigl(a=c \ \ \wedge \ \  b\leq d\bigr)\bigr].
\end{equation*}

\begin{example} \label{elo1}
Let $\cL_n$ be the subset of $\nC^n$ defined by
 \[
\cL_n := \{ (z_1, \ldots, z_n) \in \nC^n : z_1\preceq z_2 \preceq
\ldots \preceq z_n  \}.
 \]
Since $\preceq$ is a total order on $\nC$, for each $V \in \cZ_n$
the set $P^{-1}(V) \cap \cL_n$ has exactly one element. Note that
the set $\cL_n$ is not closed in $(\nC^n,\dist_{\infty})$. To show
this consider the sequence $\bigl\{\bigl(-1/k + i, 1/k - i\bigr)
\bigr\}_{k=1}^{\infty}$ in $\cL_2$ which converges to $(i,-i)$ in
$(\nC^2,\dist_{\infty})$. Clearly $(i,-i) \not\in \cL_2$. Thus
$\cL_2$ is not closed.
\end{example}

\begin{example} \label{elo2}
Define the function $L:\cZ_n\to\nC^n$ by $L :=
\bigl(P|_{\cL_n}\bigr)^{-1}$, where $\cL_n$ was defined in Example
\ref{elo1}. Thus $L(V)=(v_1,\ldots,v_n)$ where $v_1\preceq v_2
\preceq \ldots \preceq v_n$ and $V=\{v_1,\ldots,v_n\}$. We remark
that the operator $L$ is not continuous. To show this we use the
same sequence as in Example \ref{elo1} and note that
 \[
\ds\Bigl(P\bigl((-1/k + i, 1/k - i)\bigr),P\bigl((-i,i) \bigr)\Bigr)
= 1/k \to 0 \ \ \ (k \to \infty).
 \]
\end{example}

\begin{remark} \label{elo3}
With a different total order on $\nC$, for example,
\begin{equation*}
z \preceq w \ \Longleftrightarrow \  \bigl[\bigl(|z| < |w|\bigr) \
\vee \ \bigl(|z|=|w| \ \wedge \ \arg(z) \leq \arg(w)\bigr)\bigr],
\end{equation*}
the reader can create examples similar to Examples \ref{elo1} and
\ref{elo2} (with the same negative conclusions).
\end{remark}

The multiplicities of roots play an important role in the
classical statement of the continuity of roots of polynomials. The
following proposition clarifies the relation between the metric
$\ds$ and the multiplicity of the elements in a particular
multiset in $\cZ_n$.

\begin{proposition} \label{exmult}
Let  $V \in \cZ_n$  be arbitrary.  Let $v_1,\ldots, v_k$ be all
the distinct elements of \ $V$ and let $m_1, \ldots, m_k$ be their
respective multiplicities as elements of \ $V$, so that
$m_1+\cdots+m_k = n$.  Put
 \[
\eta(V) :=
 \begin{cases}
 \frac{1}{2}\min\bigl\{ |v_j - v_l|, \ j \neq l, \ j,l
\in \{1,\ldots,k \} \bigr\} & \ \ \text{for} \ \ k > 1, \\
1 & \ \  \text{for} \ \ k = 1.
\end{cases}
 \]
Then for each $U \in \cZ_n$ such that $\ds(V,U) < \eta(V)$ we have
that each disk $D(v_j,\eta(V)), \, j =1,\ldots,k$,  in the complex
plane contains exactly $m_j$ elements of $U$ counted according to
their multiplicities in $U$.
\end{proposition}
\begin{proof}
Let $U \in \cZ_n$ be such that $\ds(V,U) < \eta(V)$. Without loss
of generality, let us consider the situation around $v_1$. Let
$\sigma \in \Pi_n$ be such that $\sigma(1) = 1$ and $v_{\sigma(j)}
= v_{\sigma(1)} = v_{1}, \, j = 1,\ldots,m_1$.  By the definition
of $\ds(V,U)$, see also \eqref{eqdss}, there exists a permutation
$\tau \in \Pi_n$ such that
\[
   \max_{1\leq j \leq m_1} |v_{\sigma(1)} - u_{\tau(j)}| <
   \eta(V).
\]
Therefore all the elements $u_{\tau(j)}, \, j = 1,\ldots, m_1$, of
$U$ lie in the disk $D(v_1,\eta(V))$. Clearly, a similar statement
holds for all the other $v_j$ and since the disks $D(v_j,\eta(V)),
\, j=1,\ldots,k$, are disjoint by the definition of $\eta(V)$, the
proposition is proved.
\end{proof}

\section{Continuity} \label{sc}

In this section we prove that the function $Z : \cP_{n,1}
\rightarrow \cZ_n$ which assigns to each polynomial $p \in
\cP_{n,1}$ the multiset of its roots $Z(p) \in \cZ_n$ is a
homeomorphism between the corresponding metric spaces.

The next theorem is the classical Cauchy inequality.  Cauchy's
result is restated in terms of the metrics introduced above to
emphasize its topological meaning.  We reproduce the simple proof
of this fact as it is found in Marden's book \cite[Theorem
27.2]{M}.

\begin{theorem}[Cauchy's Inequality] \label{tci}
Define $e_n \in \cP_{n,1}$ by $e_n(z):=z^n, \, z \in \nC$, and for
any $p\in\cP_{n,1}$ let $Z(p)\in\cZ_n$ be the multiset of the
roots of $p$. Then for an arbitrary polynomial $p \in \cP_{n,1}$
we have
\begin{equation}
\ds\bigl(O,Z(p)\bigr) < 1 + \dpo(e_n,p).
\end{equation}
(Recall that by $O\in\cZ_n$ we denote the multiset of $n$ zeros.)
\end{theorem}

\begin{proof}
Let $p(z)=z^n+a_{n-1}z^{n-1}+\ldots+a_1z+a_0 \in \cP_{n,1}$ and
let $Z(p) = \{ z_1,\ldots,z_n\}$ be the roots of $p$. The theorem
claims that the following inequality holds:
\begin{equation} \label{cauchy}
    \max_{1\leq j\leq n} |z_j| <
                    1 + \max_{0\leq j\leq n-1} |a_j|\,.
\end{equation}
Let $c := \max\{|a_j|:0\leq j\leq n-1\} = \dpo(e_n,p)$. First
notice that if any root $|z_k|\leq 1$ then the inequality $|z_k| <
1 + \max\{|a_j|:0\leq j\leq n-1\}$ is trivially satisfied. Now let
$z\in \nC, \, |z|>1$. We have
\begin{eqnarray*}
|p(z)| &\geq& |z|^n - \sum_{j=0}^{n-1} |a_j||z|^j\\
       &\geq& |z|^n\Bigl(1 - c \sum_{j=1}^n |z|^{-j}\Bigr)\\
       &>& |z|^n\Bigl(1 - c \sum_{j=1}^\infty |z|^{-j}\Bigr)\\
       &>& |z|^n\Bigl(1-\frac{c}{|z|-1}\Bigr) =
       |z|^n\frac{|z|-(1+c)}{|z|-1}.
\end{eqnarray*}
Therefore, if we actually have $|z|> 1+c$, then $|p(z)|>0$ and $z$
cannot be one of the roots of $p$. This means that all roots of
$p$ must satisfy inequality \eqref{cauchy}.
\end{proof}
As an immediate consequence we have:
\begin{corollary} \label{cRbdd}
The function $Z: \cP_{n,1} \rightarrow \cZ_n$ maps each bounded
set in $(\cP_{n,1}, \dpo)$ into a bounded set in $(\cZ_n,
\ds)$.\qed
\end{corollary}

As we did above, by $\Pi_n$ we denote the set of all permutations
of $\{1,\ldots,n\}$. In the following theorem and in Section
\ref{srCn} we shall use the notation:
 \begin{equation} \label{eqdsig}
\mb{u}_{\sigma} := \bigl(u_{\sigma(1)},\ldots,u_{\sigma(n)}\bigr),
 \ \ \ \text{for} \ \ \ \sigma \in \Pi_n, \ \ \
        \mb{u} =(u_1,\ldots,u_n) \in \nC^n.
 \end{equation}

\begin{theorem} \label{tvpc}
The function $\Phi: \cZ_n \to \cP_{n,1}$ defined by
\[
 \Phi\bigl(\{z_1,\ldots,z_n\}\bigr) := \prod_{j=1}^n (z-z_j)\,,
\]
is a continuous function between $(\cZ_n, \ds)$ and $(\cP_{n,1},
\dpo)$.
\end{theorem}
\begin{proof}
Let $\{z_1,\ldots,z_n\}\in\cZ_n$ be the roots of $p(z)=z^n
+a_{n-1}z^{n-1} +\cdots+a_1z+a_0 \in \cP_{n,1}$. By Vi\`ete's
formulas,
\begin{alignat*}{2}
   a_0 &= \ (-1)^nz_1z_2\cdots z_n & \ =:& \ \psi_1(z_1,\ldots,z_n)\\
   a_1 &= \ (-1)^{n-1}\sum_{k=1}^n\prod_{j\not=k} z_j
   & \ =:& \ \psi_2(z_1,\ldots,z_n)\\
   \vdots\  &=  \hspace{20mm} \vdots & \ =:&
             \hspace{15mm} \vdots \hspace{1mm} \\
   a_{n-1} &= \ -(z_1+z_2+\ldots+z_n) & \ =:& \
   \psi_n(z_1,\ldots,z_n).
\end{alignat*}
As a linear combination of products of continuous functions, each
function $\psi_k:\nC^n \to \nC, \, k=1,\ldots,n$, is continuous.
Also note that each function $\psi_k$ is symmetric, that is
\begin{equation*}
\psi_k(\mb{u}) = \psi_k(\mb{u}_{\sigma}), \ \ \ \text{for all} \ \
\ \mb{u} \in \nC^n, \ \ \sigma \in  \Pi_n, \ \ k \in
\{1,\ldots,n\}.
\end{equation*}
(In fact each $\psi_k$ is a constant multiple of an elementary
symmetric polynomial.)

Consider the function $\Psi:\nC^n \to \nC^n$ defined by
\begin{equation*}
\Psi(\mb{v}) = \bigl(\psi_1(\mb{v}),
 \ldots,\psi_n(\mb{v})  \bigr), \ \ \ \ \mb{v} \in \nC^n.
\end{equation*}
The function $\Psi:\nC^n \to \nC^n$ is continuous and symmetric,
since each of its components $\psi_k$ is continuous and symmetric.
Therefore for each $\epsilon > 0$ and each $\mb{v} \in \nC^n$
there exists $\delta(\epsilon,\mb{v}) > 0$ such that
\begin{equation*} 
\mb{w} \in \nC^n, \ \ \dist_{\infty}(\mb{v},\mb{w}) <
\delta(\epsilon,\mb{v}) \ \ \ \Longrightarrow \ \ \
\dist_{\infty}\bigl(\Psi(\mb{v}),\Psi(\mb{w})\bigr) < \epsilon.
\end{equation*}
Also
\begin{equation*} 
\Psi(\mb{u}) = \Psi(\mb{u}_{\sigma}), \ \ \ \text{for all} \ \ \
\mb{u} \in \nC^n, \ \ \sigma \in  \Pi_n.
\end{equation*}
The last two displayed relations yield
\begin{equation} \label{eqPsics}
\mb{w} \in \nC^n, \
\min_{\sigma\in\Pi_n}\dist_{\infty}(\mb{v},\mb{w}_{\sigma}) <
\delta(\epsilon,\mb{v}) \ \ \Longrightarrow \ \
\dist_{\infty}\bigl(\Psi(\mb{v}),\Psi(\mb{w})\bigr) < \epsilon.
\end{equation}

Let $K:\cZ_n \to \nC^n$ be the function defined in \eqref{eqdK}
and let $V, W \in \cZ_n$. By the definition of $\ds$ and
\eqref{eqdsig} we have
\begin{equation*}
\ds(V,W) = \min_{\sigma\in\Pi_n}
 \dist_{\infty}\bigl(K(V),K(W)_{\sigma}\bigr).
\end{equation*}
With this observation, \eqref{eqPsics} yields
\begin{equation} \label{eqPhi1}
\begin{split}
W \in \cZ_n, \ \ \ds(V,W) < & \,\,\delta(\epsilon,K(V)) \\ & \ \ \
\Longrightarrow \ \ \dist_{\infty}\bigl(\Psi(K(V)),\Psi(K(W))\bigr)
< \epsilon.
\end{split}
\end{equation}
The definitions of $\Phi$ and $\Psi$ and the proof of Proposition
\ref{ePn} imply that
\begin{equation} \label{eqPhi2}
\dist_{\infty}\bigl(\Psi(K(V)),\Psi(K(W))\bigr) = \dpo
\bigl(\Phi(V),\Phi(W)\bigr), \ \ \ V,W \in \cZ_n.
\end{equation}
Substituting \eqref{eqPhi2} in \eqref{eqPhi1} we get that for each
$\epsilon >0$ and each $V\in\cZ_n$ there exists
$\delta(\epsilon,K(V)) > 0$ such that
\begin{equation*}
W \in \cZ_n, \ \ds(V,W) < \delta(\epsilon,K(V)) \  \Longrightarrow \
\dpo \bigl(\Phi(V),\Phi(W)\bigr) < \epsilon.
\end{equation*}
This proves the continuity of $\Phi$.
\end{proof}

Now we can prove that the space of roots and the space of
polynomials are homeomorphic.

\begin{theorem} \label{tmt}
The function $Z : \cP_{n,1} \rightarrow \cZ_n$ which associates
with each polynomial $p \in \cP_{n,1}$ the multiset of its roots
$Z(p) \in \cZ_n$ is a homeomorphism between $(\cP_{n,1}, \dpo)$
and $(\cZ_n, \ds)$.
\end{theorem}
\begin{proof}
Clearly the functions $Z$ and $\Phi$ are each other's inverse, and
so $\Phi: \cZ_n \rightarrow \cP_{n,1}$ is a bijection. Let us
verify the assumptions of Theorem \ref{tc}:

\noindent (a) \ By Theorem \ref{tbscs}, each bounded and closed
subset of the metric space $\bigl( \cZ_n, \ds\bigr)$ is compact.

\noindent (b) \ By Theorem \ref{tvpc}, $\Phi$ is continuous.

\noindent (c) \ By Corollary \ref{cRbdd}, the function $\Phi^{-1}
= Z$ maps bounded subsets of $\cP_{n,1}$ into bounded subsets of
$\cZ_n$.

Thus Theorem \ref{tc} applies and we conclude that $\Phi^{-1} = Z$
is continuous. Consequently $Z$ is homeomorphism and theorem is
proved.
\end{proof}

\section{Roots in $\nC^n$} \label{srCn}

In Section \ref{smsor} we introduced a bijection $K$ between
$\cZ_n$ and a subset $\cK$ of $\nC^n$ such that for each $V \in
\cZ_n$ the $n$-tuple $K(V)$ and the multiset $V$ have the same
elements, counting multiplicities.  Example \ref{elo2} offers a
specific bijection $L$ between $\cZ_n$ and a subset $\cL_n$ of
$\nC^n$. This bijection turns out not to be continuous. Since the
space $\nC^n$ is more familiar than $\cZ_n$, it would be desirable
to have a bijection $K:\cZ_n \to \cK \subset \nC^n$ which is a
homomorphism between $(\cZ_n,\ds)$ and $(\cK,\dist_{\infty})$. In
this section we prove that this is not possible.

\begin{theorem} \label{nocont}
Let $P$ be defined by \eqref{eqdP}. Let $\cK$ be a subset of
$\nC^n$ with the property that for each $V \in \cZ_n$ the set $\cK
\cap P^{-1}(V)$ has exactly one element.  Let $K:\cZ_n \to \nC^n$
be defined by $K = \bigl(P|_{\cK}\bigr)^{-1}$.  Then $K$ is
continuous if and only if its range $\cK$ is closed in
$(\nC^n,d_\infty)$.
\end{theorem}

\begin{proof}
Assume that $K$ is continuous. Let $\{\mathbf u_k\}$ be a Cauchy
sequence in $\cK$. Since the function $P$ satisfies
\eqref{dsineq}, the sequence $\{P(\mathbf u_k)\}$ is Cauchy in
$\cZ_n$. As $\cZ_n$ is complete by Theorem \ref{tbscs} and
Proposition \ref{p13}, this sequence is convergent, say, to $V$ in
$(\cZ_n,\ds)$.  Since $K$ is continuous the sequence $\{\mathbf
u_k\}=\{K(P(\mathbf u_k))\}$ converges to $K(V) \in \cK$. Thus
$\cK$ is closed in $(\nC,d_\infty)$.

To prove the converse assume that $\cK$ is closed. Then the
function $P|_{\cK}:\cK \to \cZ_n$ satisfies all the assumptions of
Theorem \ref{tc} (recall that $\cK$ is equipped with the metric
$\dist_{\infty}$).  Assumption (\ref{tca}) in Theorem \ref{tc} is
satisfied since each bounded and closed subset of $\cK$ is bounded
and closed in $(\nC^n,\dist_{\infty})$ and therefore compact in
$\nC^n$ and consequently compact in $\cK$. Assumption (\ref{tcb})
in Theorem \ref{tc} follows from \eqref{dsineq}, and (\ref{tcc})
follows from Proposition \ref{pKbb}.
\end{proof}

\begin{theorem} \label{tKnc}
Let $\cK$ be as in Theorem {\rm \ref{nocont}}. Then $\cK$ is not
closed in $(\nC^n,d_\infty)$.
\end{theorem}

\begin{proof}
Let $\cD$ be the set of all points $\mathbf
u=(u_1,\ldots,u_n)\in\nC^n$ such that $u_k\not=u_j$ whenever
$k\not=j$. For a point $\mathbf w$ in $\nC^n$ and $r > 0$ let
\[
  B(\mathbf w,r) = \bigl\{\mathbf v\in\nC^n : \dist_\infty(\mathbf
  w,\mathbf v)< r \bigr\}
\]
be the open ball centered at $\mathbf w$ and with radius $r$.
Also, define $\Pi_n^*$ to be the set of all permutations of
$\{1,\ldots,n\}$ minus the identity permutation.

By contradiction, suppose that $\cK$ is closed in
$(\nC^n,d_\infty)$. Let $\mathbf u\in\cK\cap\cD$, that is, all the
coordinates of $\mathbf u \in \cK$ are mutually distinct.  By the
definition of $\cK$, for every $\sigma\in\Pi_n^*$ we have that
${\mathbf u}_{\sigma}\in\nC^n\!\setminus\!\cK$. Since
$\nC^n\!\setminus\!\cK$ is open, there exist an $r_\sigma > 0$
such that the entire open ball $B({\mathbf u}_{\sigma},r_\sigma)$
is contained in $\nC^n\!\setminus\!\cK$. Now we put
\[
   r := \min \bigl\{r_\sigma \, : \,  \sigma \in \Pi_n^*  \bigr\}
\]
and prove that the ball $B(\mathbf u,r)$ is entirely contained in
$\cK$.  To see this, pick a $\mathbf v\in B(\mathbf u,r)$. Then by
our choice of $r$ it follows that ${\mathbf v}_{\sigma}$ is
contained in $B({\mathbf u}_{\sigma},r_{\sigma})$ (and thus
${\mathbf v}_{\sigma}\not\in\cK$) for all $\sigma\in\Pi_n^*$. Since
our construction of $\cK$ requires that some permutation of the
coordinates of $\mathbf v$ be contained in $\cK$, and the only one
we have left is $\mathbf v$ itself, we conclude that $\mathbf
v\in\cK$. So, $B(\mathbf u,r)\subset\cK$, as claimed. We have thus
proved that all the points in $\cK\cap\cD$ (i.e., those with $n$
distinct coordinates) are interior points of $\cK$.

Now let $\sigma\in\Pi_n^*$. Since $\mb{u} \in \cK \,\cap\, \cD$,
we have ${\mathbf u}_{\sigma} \in \cD \!\setminus\! \cK$. By Lemma
\ref{lDpwc}, $\cD$ is pathwise connected. Therefore there exists a
continuous function $\Theta: [0,1] \to \cD$ such that $\Theta(0) =
\mathbf u$ and $\Theta(1) = {\mathbf u}_{\sigma}$. Let
 \begin{equation} \label{eqda}
  a := \sup \, \bigl\{t \in [0,1] \, : \, \Theta(t) \in \cK \bigr\}.
 \end{equation}
This supremum exists since $\Theta(0) = \mathbf u \in \cK$ so the
set on the right-hand side of \eqref{eqda} is not empty. As we
assume that $\cK$ is closed, $\Theta(a) \in \cK$. Therefore $a <
1$. The range of $\Theta$ is a subset of $\cD$, and thus
$\Theta(a) \in \cK\cap\cD$ and consequently $\Theta(a)$ must be an
interior point of $\cK$. Since $\Theta$ is continuous this
contradicts the definition of $a$. Thus $\cK$ cannot be closed.
\end{proof}

An immediate consequence of the previous two theorems is:

\begin{corollary} \label{cnc}
The operator $K$ defined in Theorem {\rm \ref{nocont}} is not
continuous.
\end{corollary}

\begin{example} \label{elo4}
Let $\cL_n, L,$ and $Z$ be as in Examples \ref{elo1}, \ref{elo2}
and Theorem \ref{tmt}.  Then the function $L\circ Z: \cP_{n,1} \to
\cL_n \subset \nC^n$ is not continuous. For simplicity, we
consider $n=2$.  The sequence of polynomials
 \[
z^2 + 1 +2i/k -1/k^2, \ \ \  k\in\nN,
 \]
converges to $z^2 + 1$ in $\bigl(\cP_{2,1},\dpo\bigr)$, but the
sequence of lexicographically ordered pairs of their roots $(-1/k
+ i, 1/k - i), k\in \nN$, does not converge in
$(\nC^2,\dist_{\infty})$ to the pair of lexicographically ordered
roots $(-i,i)$ of $z^2+1$.
\end{example}

\begin{remark}
A metric space setting for Theorem \ref{tmt} is also provided in
\cite{NP} and parts of our proof are similar to the proofs in
\cite{NP}.  In \cite{NP} the authors consider two metric spaces: the
space of all monic polynomials of degree $n$ and the space of their
roots considered as ordered $n$-tuples of complex numbers (ordered
lexicographically as explained in Example \ref{elo1}) and equipped
with the $\dist_{\infty}$ metric.   Example \ref{elo4} points out
the difficulty with this setting (which invalidates the argument in
\cite{NP}). Moreover Corollary \ref{cnc} and Theorem \ref{tmt} imply
that it is {\em not} possible to identify the roots of monic
polynomials with unique $n$-tuples and equip such a set with the
$\dist_{\infty}$ metric and have a homeomorphism between such space
of roots and the space of polynomials.  This indicates that the
metric $\ds$ is the natural metric on the roots.
\end{remark}

\section{Final remarks} \label{sfr}

We conclude with some historical remarks. In 1939 Ostrowski
\cite{O1} published his own form of the perturbation theorem for
polynomial roots. We quote it from \cite[Appendix A]{O}.

\begin{theorem}
Consider two polynomials
\begin{align*}
f(x)& =a_0 x^n+ \cdots +a_n,   \ \ \ a_0 =1, \\
g(x) &= b_0 x^n+ \cdots + b_n, \ \ \ b_0 = 1.
\end{align*}
Let the $n$ roots of $f(x)$ be $x_1,\ldots,x_n$, those of $g(x),
\, y_1,\ldots,y_n$. Put
\[
\gamma = 2 \, \Gamma, \ \ \   \Gamma = \max_{\nu > 0} \bigl(
|a_\nu|^{1/\nu},|b_\nu|^{1/\nu} \bigr).
\]
Introduce the expression
 \[
\varepsilon = \sqrt[n]{\sum_{\nu=1}^n \, |b_\nu - a_\nu| \,
\gamma^{n-\nu}}.
 \]
 The roots $x_\nu$ and $y_\nu$ can be ordered in such a way that
 we have
  \[
|x_\nu - y_\nu| < (2n-1) \, \varepsilon \ \ \ \ (\nu =
1,\ldots,n).
  \]
\end{theorem}

We can see that Ostrowski's statement was quite ``ready'' for the
language of the metric $\ds$, as it essentially contains the
definition we give of $\ds$ in Section \ref{smsor}. To show an
alternate presentation of the classical perturbation theorem
(though this time without the kind of numerical estimate that
Ostrowski wanted to obtain), here is the one given in \cite{M}:

\begin{theorem}
Let
\begin{eqnarray*}
  f(z) &=& a_0 + a_1z + \cdots + a_nz^n =
    a_n \prod_{j=1}^p (z-z_j)^{m_j}\,,
  \quad\quad\quad
           a_n\not=0,\\
  F(z) &=& (a_0+\varepsilon_0) + (a_1+\varepsilon_1)z + \cdots +
           (a_{n-1}+\varepsilon_{n-1})z^{n-1} + a_nz^n
\end{eqnarray*}
and let
\[
  0 < r_k < \min |z_k-z_j|\,,\quad\quad\quad\quad
           j=1,2,\ldots,k-1,k+1,\ldots,p\,.
\]
Then there exists a positive number $\varepsilon$ such that, if\,
$|\varepsilon_i|\leq\varepsilon$ for $i=0,\ldots,n-1$, then $F(z)$
has precisely $m_k$ zeros in the circle $C_k$ with center $z_k$
and radius $r_k$.
\end{theorem}

As a last quote, here is a version of the continuity theorem from
the recent major survey of the theory of polynomials by Rahman and
Schmeisser \cite[Theorem 1.3.1 and Supplement]{RS}

\begin{theorem}
Let
\[
  f(z) = \sum_{\nu=0}^n a_\nu z^\nu = \prod_{j=1}^k (z-z_j)^{m_j} \ \
  \ (m_1+\cdots+m_k=n)
\]
be a monic polynomial of degree $n$ with distinct zeros
$z_1,\ldots,z_k$ of multiplicities $m_1,\ldots,m_k$. Then, given a
positive $\varepsilon<\min_{1\leq i\leq j\leq k} |z_i-z_j|/2$, there
exists a $\delta>0$ so that any monic polynomial
$g(z)=\sum_{\nu=0}^n b_\nu z^\nu$ whose coefficients satisfy
$|b_\nu-a_\nu|<\delta$, for $\nu=1,\ldots,n-1$, has exactly $m_j$
zeros in the disc
\[
  D(z_j,\varepsilon) \ \ (j=1,\ldots,k).
\]
Further, if we let
\[
   A := \max\bigl\{1, 2|a_\nu|^{1/(n-\nu)} : \nu=0,\ldots,n-1\bigr\},
\]
and let the zeros of $f$ be denoted by $\zeta_1,\ldots,\zeta_n$,
where an $m$-fold zero is now listed $m$ times, then, for
sufficiently small $\delta>0$, there exists a numbering of the zeros
of $g$ as $\omega_1,\ldots,\omega_n$ such that $\max_{1\leq\nu\leq
n} |\omega_\nu-\zeta_\nu|\leq 4A\delta^{1/n}$.
\end{theorem}

To conclude: in every case known to us, the classical perturbation
theorem has been presented as a continuity result (in a more or less
convoluted way) and it has been proved by many authors using a
variety of techniques (mostly from complex function theory, or
trying to obtain useful numerical estimates). We hope that our
topological presentation, and the emphasis on the homeomorphic
relation between roots and polynomials, may have added to the
understanding of this beautiful, age-old result.

\end{document}